\input amstex
\documentstyle{amsppt}
\def\sdp{\medspace \times \kern -1.9pt  
         \vrule width0.4pt height 4.7pt depth -0.3pt \medspace}

\magnification=\magstep1
\hsize=5.2in
\vsize=6.8in

\topmatter

\title A 1-COHOMOLOGY CHARACTERIZATION OF PROPERTY (T) \\
IN VON NEUMANN ALGEBRAS
\endtitle

\author JESSE PETERSON
\endauthor

\rightheadtext{A 1-cohomology characterization of property (T)}

\affil University of California, Los Angeles\endaffil

\address Math Dept  
UCLA, Los Angeles, CA 90095-155505\endaddress

\email jpete\@math.ucla.edu \endemail

\abstract We obtain a characterization of property (T) for von Neumann algebras
in terms of 1-cohomology similar to the Delorme-Guichardet Theorem for groups.
\endabstract



\endtopmatter

\document

\head 0. Introduction. \endhead

\vskip .1in
The analogue of group representations in von Neumann
algebras is the notion of {\it correspondences} which is due to Connes ([C2], [C3], [P1]), 
and has been a very useful in defining notions such as property (T) and amenability for 
von Neumann algebras.
It is often useful to view group
representations as positive definite functions which we obtain through a GNS construction.
Correspondences of a von Neumann algebra $N$ can also be viewed in two separate ways, as 
{\it Hilbert $N$-$N$ bimodules} $\Cal H$, or as {\it completely positive maps} 
$\phi: N \rightarrow N$, and the equivalence of these two descriptions 
is also realized via a GNS construction.  This allows one to characterize property (T) for von 
Neumann algebras in terms of completely positive maps.

For a countable group $G$ there is also a notion of conditionally negative definite functions 
$\psi:G \rightarrow \Bbb C$ which satisfy $\psi(g^{-1}) = \overline {\psi(g)}$ and the condition:  
$\forall n \in \Bbb N$, $\alpha_1,\alpha_2, \dots, \alpha_n \in \Bbb C$, $g_1,g_2, \dots,g_n \in G$,
if $\Sigma_{i=1}^n \alpha_i = 0$ then 
$\Sigma_{i,j=1}^n \overline{\alpha_j}\alpha_i\psi(g_j^{-1}g_i) \leq 0$.
Real valued conditionally negative definite functions can be viewed as cocycles 
$b \in B^1(G,\pi)$ where $\pi:G \rightarrow \Cal O(\Cal H)$ is an orthogonal representation
of $G$ (see [BdHV]).  Real valued conditionally negative definite functions can also be viewed as
generators of semigroups of positive definite functions by Schoenberg's Theorem.  These equivalences then make
it possible for certain connections between 1-cohomology,
conditionally negative definite functions, and positive definite deformations,
for example the Delorme-Guichardet Theorem [De] [G] which states that a 
group has property (T) of Kazhdan [Ka] if and only if the first cohomology vanishes for 
any unitary representation.

It was Evans who introduced the notion of bounded conditionally completely positive/negative maps [E] 
related to the
study the infinitesimal generators of norm continuous semigroups of completely positive maps.
He noted that this definition gives an analogue to conditionally positive/negative definite functions on groups.
We will extend the notion of conditionally completely negative maps to unbounded maps 
and use a GNS type construction to alternately view them as closable derivations into a Hilbert $N$-$N$ bimodule.
This is done in the same spirit as ([S1],[S2]) where Sauvageot makes a connection between
quantum Dirichlet forms, and differential calculus.  Indeed, it is shown in 1.7.1 that
conditionally completely negative maps are in fact extensions of generators associated
to completely Dirichlet forms, however we are coming from a different 
perspective here and so we will present the correspondence between conditionally completely negative maps and
closable derivations in a way more closely related to group theory.

In studying various properties of groups such as property (T) or the Haagerup property one can 
give a characterization of these properties in terms of boundedness conditions on conditionally negative definite
functions (e.g. [AW]), hence one would hope that this is possible for von Neumann algebras as well.

We will show that one can indeed obtain a characterization
of property (T) in this way.
The main result is that a separable finite factor has property (T) 
if and only if the 1-cohomology spaces of closable derivations vanish
whenever the domain contains a non-$\Gamma$ set (see sec. 3 for the definition 
of a non-$\Gamma$ set).  

\proclaim{0.1. Theorem}Suppose that $N$ is a separable finite factor, then 
the following conditions are equivalent:

\vskip .05in
\noindent
(1)  $N$ has property (T).

\vskip .05in
\noindent
(2)  $N$ does not have property $\Gamma$ and given any weakly dense $*$-subalgebra $N_0 \subset N$, $1 \in N_0$ such that
$N_0$ contains a non-$\Gamma$ set, we have that
every densely defined closable derivation on $N_0$ into a Hilbert $N$-$N$ bimodule is inner.

\vskip .05in
\noindent
(3)  There exists a weakly dense $*$-subalgebra $N_0 \subset N$ such that $N_0$ is countably generated
as a vector space and every closable derivation into a Hilbert $N$-$N$ bimodule whose domain contains $N_0$ is inner.
\endproclaim

This is the analogue to the Delorme-Guichardet Theorem for groups.  As a corollary we obtain that if 
$X_1, \ldots, X_n$ generate a finite factor with property (T), 
and if at least one of the $X_j$'s has diffuse spectrum then the derivations $\partial_{X_i}$ from [V1] cannot all be 
closable and hence the conjugate variables cannot all exist in $L^2(N, \tau)$.

\proclaim{0.2. Corollary}  Suppose that $N$ is a finite factor with property (T), let 
$X_1, \ldots, X_n$ generate $N$ as a von Neumann algebra
such that $\Bbb C[X_1, \ldots, X_n]$ contains a non-$\Gamma$ set.
If one of the $X_j$'s has diffuse spectrum then 
$\Phi^*(X_1, \ldots, X_n) = \infty$.
\endproclaim

We also give an application 
showing that many amalgamated free products of finite von Neumann algebras do not have property (T).

\proclaim{0.3. Theorem} Let $N_1$ and $N_2$ be finite von Neumann algebras with 
with normal faithful tracial states $\tau_1$ and $\tau_2$ respectively, suppose that $B$ is a common
von Neumann subalgebra such that $\tau_1|_B = \tau_2|_B$.  If there are unitaries
$u_i \in \Cal U (N_i)$ such that $E_B(u_i) = 0$, $i = 1,2$.
Then $M = N_1 *_B N_2$ does not have property (T).
\endproclaim

Other than the introduction there are 4 sections.  Section 1 establishes the definitions and notations as well as 
gives the connection between closable derivations, conditionally completely negative maps, and semigroups of 
completely positive maps.  In section 2 we characterize when a closable derivation is inner in terms of the 
conditionally completely negative map and the semigroup.  In Section 3 we state and prove the main theorem 
(3.2), and in section 4 we give the application with amalgamated free products (4.2).

\vskip .1in
\head 1. A GNS-type construction. \endhead

\vskip .1in
\noindent
{\bf 1.1. Conditionally completely negative maps}. Let $N$ be a finite von Neumann algebra 
with normal faithful trace $\tau$.

\proclaim{Definition} Suppose $\Psi:N \rightarrow L^1(N, \tau)$ is a $*$-preserving
linear map whose domain is a weakly dense $*$-subalgebra $D_\Psi$ of $N$ such that
$1 \in D_\Psi$, then $\Psi$ is a conditionally completely negative (c.c.n.) map on $N$
if the following condition is satisfied:

\vskip .05in
\noindent  
$(1.1.1)$.  $\forall n \in \Bbb N$, $x_j, y_j \in D_\Psi$, $j \leq n$,
if $\Sigma_{j=1}^n x_jy_j = 0$ then $\Sigma_{i,j=1}^n y_j^*\Psi(x_j^*x_i)y_i \leq 0$. 
  
\endproclaim

It is not hard to see that condition $(2.1.1)$ can be replaced with the condition:

\vskip .05in
\noindent 
$(1.1.1)'$. $\forall n \in \Bbb N$, $x_j, y_j \in D_\Psi$, $j \leq n$,
if $\Sigma_{j=1}^n x_jy_j = 0$ then $\Sigma_{i,j=1}^n \tau(\Psi(x_j^*x_i)y_iy_j^*) \leq 0$. 

\vskip .1in
If $\phi:N \rightarrow N$ is a completely positive map and $k \in N$ then 
$\Psi(x) = k^*x + xk - \phi(x)$ gives a 
map which is c.c.n. and bounded.  If 
$\delta:N \rightarrow L^2(N, \tau)$ is a derivation then $\delta$ is c.c.n.  Also if 
$\Psi$ is a c.c.n. map and $\alpha:N \rightarrow N$ is a $\tau$-preserving automorphism then
$\Psi' = \alpha \circ \Psi \circ \alpha^{-1}$ is another c.c.n. map.

One can check that if $\Psi_1$ and $\Psi_2$ are 
c.c.n. such that $D_{\Psi_1} \cap D_{\Psi_2}$ is weakly dense in $N$, and if $s,t \geq 0$, then 
$\Psi = s\Psi_1 + t\Psi_2$ is c.c.n.  Also if $\{\Psi_t\}_t$ is a 
family of c.c.n. maps on the same domain and $\Psi$ is the pointwise $\| \cdot \|_1$-limit of $\{\Psi_t\}_t$ then
$\Psi$ is c.c.n.

We say that $\Psi$ is symmetric if $\tau(\Psi(x)y) = \tau(x\Psi(y)), \forall x,y \in D_\Psi$.
We say that $\Psi$ is conservative if $\tau \circ \Psi = 0$.  We also say that $\Psi$ is closable
if the quadratic form $q$ on $L^2(N, \tau)$ given by $D(q) = D_\Psi$, $q(x) = \tau(\Psi(x)x^*)$ 
is closable.  Note that we will see in (1.3) that if $\Psi: D_\Psi \rightarrow L^2(N, \tau) \subset L^1(N, \tau)$
is a conservative symmetric c.c.n. map then $\Psi$ is automatically closable.

Note that if $\Psi$ is a conservative symmetric c.c.n. map then 
$\tau(\Psi(1)x) = \tau(\Psi(x)) = 0, \forall x \in D_\Psi$, hence $\Psi(1) = 0$.  
Also note that if $\Psi$ is symmetric and 
$\Psi(1) \geq 0$ then given any $x \in D_\Psi$, if we
let $x_1 = x, x_2 = 1, y_1 = -1, y_2 = x$, then the above condition implies that 
$\tau(\Psi(x)x^*) \geq 0$, so that we actually have positivity instead of just the 
symmetry condition.

\vskip .1in
\noindent
{\bf 1.2. Closable derivations}.  Let $\Cal H$ be a Hilbert $N$-$N$ bimodule, a derivation of $N$ is a 
(possibly unbounded) map $\delta: N \rightarrow \Cal H$
which is defined on a weakly dense $*$-subalgebra $D_\delta$ of $N$ such that $1 \in D_\delta$,
and such that $\forall x, y \in D_\delta$, $\delta(xy) = x\delta(y) + \delta(x)y$.
$\delta$ is closable if it is closable as an operator from $L^2(N, \tau)$ to $\Cal H$.

$\delta$ is inner if $\delta(x) = x\xi - \xi x$ for some $\xi \in \Cal H$.
$\delta$ is spanning if $\overline{sp}D_\delta\delta(D_\delta) = \Cal H$. 
$\delta$ is real if $\langle x\delta(y),\delta(z) \rangle_{\Cal H} = 
\langle \delta(z^*), \delta(y^*)x^* \rangle_{\Cal H}$, $\forall  x,y,z \in D_\delta$.

If $\delta':D_\delta \rightarrow \Cal H'$ is another derivation then we say that
$\delta$ and $\delta'$ are equivalent if there exists
a unitary map $U:\Cal H \rightarrow \Cal H'$ such that $U(x\delta(y)z) = xU(\delta(y))z = x\delta'(y)z$
for all $x,y,z \in D_\delta$. 

Recall that if $\Cal H$ is a Hilbert $N$-$N$ bimodule then we can define the adjoint bimodule 
$\Cal H^o$ where $\Cal H^o$ is the conjugate Hilbert space of $\Cal H$ and the bimodule
structure is given by $x \xi^o y = (y^* \xi x^*)^o$.
If $\delta:D_\delta \rightarrow \Cal H$ is a closable derivation then we may define the 
adjoint derivation $\delta^o:D_\delta \rightarrow \Cal H^o$ by setting
$\delta^o(x) = \delta(x^*)^o$, then $\delta^o$ is a closable derivation
and furthermore the derivations ${1 \over 2} (\delta + \delta^o)$, and
${1 \over 2} (\delta - \delta^o)$ are real derivations from $D_\delta$ to $\Cal H \oplus \Cal H^o$.

\vskip .1in
\noindent
{\bf 1.3. From conditionally completely negative maps to closable derivations}.
Let $\Psi$ be a conservative symmetric c.c.n. map on $N$ with domain $D_\Psi$.
We associate to $\Psi$ a derivation in the following way (compare with [S1]): 

let $\Cal H_0 = \{\Sigma_{i=1}^n x_i\otimes y_i \in D_\Psi\otimes D_\Psi| \Sigma_{i=1}^n x_iy_i = 0 \}$.
Define a sesquilinear form on 
$\Cal H_0$ by $\langle \Sigma_{i=1}^n x_i'\otimes y_i', \Sigma_{j=1}^m x_j\otimes y_j \rangle_\Psi = 
-{1 \over 2}\Sigma_{i=1}^n \Sigma_{j=1}^m \tau(\Psi(x_j^*x_i')y_i'y_j^*)$.  The positivity of 
$\langle \cdot,\cdot \rangle_\Psi$ is equivalent to the c.c.n. condition on $\Psi$.
Let $\Cal H$ be the closure of $\Cal H_0$ after we mod out by the kernel of $\langle \cdot,\cdot \rangle_\Psi$.
If $p = \Sigma_{k=1}^n x_k \otimes y_k$ such that $\Sigma_{k=1}^n x_ky_k = 0$ then 
$x \mapsto -{1 \over 2}\Sigma_{i,j=1}^n \tau(x_j^*xx_i\Psi(y_iy_j^*))$ is a positive normal functional
on $N$ with norm $\langle p,p \rangle_\Psi$.  Similarly 
$y \mapsto -{1 \over 2}\Sigma_{i,j=1}^n \tau(\Psi(x_j^*x_i)y_iyy_j^*)$ is a positive normal functional
on $N$ with norm $\langle p,p \rangle_\Psi$.  We also have left and right commuting actions of $D_\Psi$
on $\Cal H_0$ given by $xpy = x(\Sigma_{k=1}^n x_k \otimes y_k)y = \Sigma_{k=1}^n (xx_k) \otimes (y_ky)$,
and by the preceeding remarks we have 
$\langle xp,xp \rangle_\Psi = \langle x^*xp,p \rangle_\Psi \leq 
  \|x^*x\|\langle p,p \rangle_\Psi = \|x\|^2\langle p,p \rangle_\Psi$ and
$\langle py,py \rangle_\Psi \leq \|y\|^2\langle p,p \rangle_\Psi$ for all $x,y \in D_\Psi$.
Hence the above actions of $D_\Psi$ pass to commuting left and right actions on $\Cal H$, 
and they extend to left and right
actions of $N$ on $\Cal H$ given by the formulas:
$$
\langle x[\Sigma_{i=1}^n x_i' \otimes y_i'],[\Sigma_{j=1}^m x_j \otimes y_j] \rangle_{\Cal H} =
\Sigma_{i=1}^n \Sigma_{j=1}^m \tau(x_j^*xx_i'\Psi(y_i'y_j^*)),
$$
$$
\langle [\Sigma_{i=1}^n x_i' \otimes y_i']y,[\Sigma_{j=1}^m x_j \otimes y_j] \rangle_{\Cal H} =
\Sigma_{i=1}^n \Sigma_{j=1}^m \tau(\Psi(x_j^*x_i')y_i'yy_j^*).
$$
Since the above forms are normal the left and right actions commute and are normal 
thus making $\Cal H$ into a Hilbert $N$-$N$ bimodule.

Define $\delta_\Psi:D_\Psi \rightarrow \Cal H$ to be given by
$\delta_\Psi(x) = [x\otimes 1 - 1\otimes x]$.
Then $\delta_\Psi$ is a derivation such that
$$
\langle \delta_\Psi(x),\delta_\Psi(y) \rangle_{\Cal H} = \langle x\otimes 1 - 1\otimes x, y\otimes 1 - 1 \otimes y \rangle_\Psi 
$$
$$
= -{1 \over 2}\tau(\Psi(y^*x)) + {1 \over 2}\tau(\Psi(x)y^*) + {1 \over 2}\tau(\Psi(y^*)x) - {1 \over 2}\tau(\Psi(1)xy^*)
$$
$$
= \tau(\Psi(x)y^*),
$$
for all $x,y \in D_\Psi$.
Also $\delta_\Psi$ is real since
$$
\langle x\delta_\Psi(y), \delta_\Psi(z) \rangle_{\Cal H} =
\langle xy \otimes 1 - x \otimes y, z \otimes 1 - 1 \otimes z \rangle_\Psi
$$
$$
= -{1 \over 2}\tau(\Psi(z^*xy)) + {1 \over 2}\tau(\Psi(xy)z^*) + {1 \over 2}\tau(\Psi(z^*x)y) - {1 \over 2}\tau(\Psi(x)yz^*) 
$$
$$
= -{1 \over 2}\tau(\Psi(1)z^*xy) + {1 \over 2}\tau(\Psi(z^*)xy) + {1 \over 2}\tau(\Psi(y)z^*x) - {1 \over 2}\tau(\Psi(yz^*)x)
$$
$$
= \langle 1 \otimes z^* - z^* \otimes 1, 1 \otimes y^*x^* - y^* \otimes x^* \rangle_\Psi =
\langle \delta_\Psi(z^*), \delta_\Psi(y^*)x^* \rangle_{\Cal H},
$$
for all $x,y,z \in D_\Psi$.
If $\Psi$ is closable then it follows that $\delta_\Psi$ is closable. 
Also note that if $\Psi : D_\Psi \rightarrow L^2(N, \tau) \subset L^1(N, \tau)$ then we would have
that $D_\Psi = D(\delta_{\Psi}^* \delta_{\Psi})$ which would show that $\delta$ (and hence also $\Psi$)
is closable.

Note that we will assume in addition that $\delta_\Psi$ is spanning by restricting ourselves to 
$\overline{sp}D_\Psi\delta(D_\Psi) \subset \Cal H$.

Also note that the requirement that $\Psi(1) = 0$ is not really much of a restriction since 
if $\Psi$ is any symmetric c.c.n. map with $\Psi(1) \in L^2(N, \tau)$
then $\Psi'(x) = \Psi(x) - {1 \over 2}\Psi(1)x - {1 \over 2}x\Psi(1)$ defines
a symmetric c.c.n. map with $\Psi'(1) = 0$.
\vskip .1in
\noindent
{\bf 1.4. From closable derivations to conditionally completely negative maps}.
Let $\Cal H$ be a Hilbert $N$-$N$ bimodule and suppose that $\delta: N \rightarrow \Cal H$
is a closable real derivation defined on a weakly dense $*$-subalgebra $D_\delta$ of $N$
with $1 \in D_\delta$.

Let $D_{\Psi} = \{ x \in D(\overline {\delta}) \cap N | y \mapsto \langle \overline{\delta}(x), \delta(y^*) \rangle 
			{\text { gives a normal linear functional on }} N \}$
then by [S2] and [DL] $D(\overline {\delta}) \cap N$ is a $*$-subalgebra and hence 
one can show that $D_{\Psi}$ is a $*$-subalgebra of $N$.
We define the map $\Psi_\delta: D_{\Psi} \rightarrow L^1(N, \tau)$
by letting $\Psi_\delta(x)$ be the Radon-Nikodym derivative of the normal linear functional
$y \mapsto \langle \overline{\delta}(x), \delta(y^*) \rangle$.
Since $\delta$ is closable $\Psi_\delta$ is also closable.

As $\delta$ is real $\Psi_\delta$ is a symmetric $*$-preserving map such that $\tau \circ \Psi = 0$ and if 
$n \in \Bbb N$, $x_1,x_2,\dots,x_n,y_1,y_x,\dots,y_n \in D_\Psi,$ such that 
$\Sigma_{i=1}^n x_iy_i = 0$ then:
$$
\Sigma_{i,j=1}^n \tau(\Psi(x_j^*x_i)y_iy_j^*) = 
 \Sigma_{i,j=1}^n \langle \delta(x_j^*x_i),\delta(y_jy_i^*) \rangle_{\Cal H}
$$
$$
= \Sigma_{i,j=1}^n \langle x_j^*\delta(x_i),y_j\delta(y_i^*) + \delta(y_j)y_i^* \rangle_{\Cal H} +
                 \langle \delta(x_j^*)x_i,y_j\delta(y_i^*) + \delta(y_j)y_i^* \rangle_{\Cal H}
$$
$$
= \Sigma_{i,j=1}^n \langle \delta(x_i)y_i,x_j\delta(y_j) \rangle_{\Cal H} +
                  \langle x_i\delta(y_i),\delta(x_j)y_j \rangle_{\Cal H}
$$
$$
= -2\| \Sigma_{i=1}^n \delta(x_i)y_i \|_{\Cal H}^2 \leq 0.
$$

Hence $\Psi_\delta$ is a conservative symmetric c.c.n. map on $D_\Psi$.

Note that if we restrict ourselves to closable derivations which are spanning then an easy calculation shows that the 
constructions above are inverses of each other in the sense that $\Psi_{\delta_\Psi}|_{D_\Psi} = \Psi$ and 
$\overline {\delta_{\Psi_\delta}} \cong \overline {\delta}$.

\vskip .1in
\noindent
{\bf 1.5. Closable derivations and c.c.n. maps from groups}.
Let $\Gamma$ be a discrete group, $(C,\tau_0)$ a finite von Neumann algebra with a normal faithful trace, and 
$\sigma$ a cocycle action of $\Gamma$ on $(C,\tau_0)$ by $\tau_0$-preserving automorphisms.  Denote by 
$N = C \sdp_\sigma \Gamma$ the corresponding cross-product algebra with trace $\tau$ given by 
$\tau( \Sigma c_g u_g ) = \tau_0(c_e)$, where $c_g \in C$ and $\{u_g\}_g \subset N$ denote the 
canonical unitaries implementing the action $\sigma$ on $C$.

Let $(\pi_0,\Cal H_0)$ be a unitary or orthogonal representation of $\Gamma$, and let 
$b:\Gamma \rightarrow \Cal H_0$ be an (additive) cocycle of $\Gamma$, i.e. 
$b(gh) = \pi_0(g)b(h) + b(g)$, $\forall g,h \in \Gamma$.  Set $\Cal H_{\pi_0}$ to be the
Hilbert space $\Cal H_0 \overline \otimes_{\Bbb R} L^2(N,\tau)$ if $\pi_0$ is an orthogonal representation and
$\Cal H_0 \overline \otimes_{\Bbb C} L^2(N,\tau)$ if $\pi_0$ is a unitary representation.  
We let $N$ act on the right of $\Cal H_{\pi_0}$
by $(\xi \otimes \hat x)y = \xi \otimes (\hat {xy})$, $x,y \in N, \xi \in \Cal H_0$ and on the left by
$c(\xi \otimes \hat x) = \xi \otimes (\hat {cx}), u_g(\xi \otimes \hat x) = (\pi_0(g)\xi) \otimes (\hat {u_gx})$,
$c \in C, x \in N, g \in \Gamma, \xi \in \Cal H_0$.  
Let $D_\Gamma$ be the $*$-subalgebra generated by $C$ and $\{u_g\}_g$, we 
define $\delta_b$ by $\delta_b(c_gu_g) = c_g\delta_b(u_g) = b(g) \otimes \hat{c_gu_g}$, $c_g \in C, g \in \Gamma$,
then we can extend $\delta_b$ linearly so that $\delta_b$ is a derivation on $D_\Gamma$.  
If $(\pi_0,\Cal H_0)$ is an orthogonal representation and
$1_g$ denotes the Dirac delta function at $g$ then:
$$
\langle cu_g\delta_b(u_h),\delta_b(u_k) \rangle = 
 \langle \pi_0(g)b(h), b(k) \rangle \langle \hat{cu_gu_h}, \hat u_k \rangle
$$
$$
= \langle -\pi_0(g)\pi_0(h)b(h^{-1}),-\pi_0(k)b(k^{-1}) \rangle \langle \hat{cu_gu_h}, \hat u_k \rangle 1_k(gh)
$$
$$
= \langle b(k^{-1}),b(h^{-1}) \rangle \langle \hat u_k^*, \hat{u_h^*u_g^*c^*} \rangle
$$
$$
= \langle \delta_b(u_k^*),\delta_b(u_h^*)u_g^*c^* \rangle,
$$
for all $g,h,k \in \Gamma$, $c \in C$, thus showing that $\delta_b$ is real.

Also we have:
$$
|\langle \delta_b(c_gu_g),\delta_b(\Sigma_{h \in \Gamma}d_hu_h) \rangle | = 
|\Sigma_{h \in \Gamma} \langle b(g),b(h) \rangle \langle \hat {c_gu_g}, \hat {d_hu_h} \rangle |
$$
$$
= \|b(g)\|^2 |\langle \hat{c_gu_g}, \Sigma_{h \in \Gamma} \hat{d_hu_h} \rangle | \leq
\|b(g)\|^2 \|c_g\| \| \Sigma_{h \in \Gamma} d_hu_h\|_1,
$$
for all $g \in \Gamma$, $c_g \in C$, $\Sigma_{h \in \Gamma} d_hu_h \in D_\Gamma$.
Hence if $x = \Sigma_{g \in \Gamma}c_gu_g \in  D_\Gamma$, and $y \in D_\Gamma$ then
$|\langle \delta_b(x), \delta_b(y) \rangle| \leq (\Sigma_{g \in \Gamma} \|b(g)\|^2 \|c_g\|) \|y\|_1$.
In particular this shows that $\delta_b$ is closable.

Now suppose that $\psi:\Gamma \rightarrow \Bbb C$ is a real valued conditionally negative definite function on 
$\Gamma$ such that $\psi(e) = 0$
and let $(\pi_\psi,b_\pi)$ be the representation and cocycle which correspond to $\psi$ through the GNS 
construction [BdHV].  Let $(\Cal H,\delta)$ denote the Hilbert $N$-$N$ bimodule and closable derivation constructed out
of $(\pi_\psi,b_\pi)$ as above and let $\Psi$ be the symmetric c.c.n. map associated to 
$(\Cal H,\delta)$ as in 1.4.  Then a calculation shows that 
$\Psi(\Sigma_gc_gu_g) = \Sigma_g \psi(g)c_gu_g$, and in fact it is an easy exercise to show that even if $\psi$ is
not real valued $\Psi(\Sigma_gc_gu_g) = \Sigma_g \psi(g)c_gu_g$ still describes a c.c.n. map.

Conversely, if $(\Cal H,\delta)$ is a Hilbert $N$-$N$ bimodule and a closable derivation such that 
$\delta$ is defined on the $*$-subalgebra generated by $C$ and $\{u_g\}_g$ then we can associate to
it a representation $\pi_0$ on $\Cal H_0 = \overline{sp}\{\delta (u_g)u_g^*|g \in \Gamma\}$ by
$\pi_0(g)\xi' = u_g\xi'u_g^*$, $\xi' \in \Cal H_0$.  Also we may associate to $\delta$ a group cocycle $b$ on 
$\Gamma$ by $b(g) = \delta(u_g)u_g^*$, $g \in \Gamma$.  If $\Psi$ is a c.c.n. map which is
also defined on the $*$-subalgebra generated by $C$ and $\{u_g\}_g$ then we can associate to it
a conditionally negative definite function $\psi$ by $\psi(g) = \tau(\Psi(u_g)u_g^*)$.  Furthermore 
if $\delta$ is real then by taking only the real span above we have that $\Cal H_0$ is a real 
Hilbert space and $\pi_0$ is an orthogonal representation, also $\psi$ is real valued if
and only if $\Psi$ is symmetric, and if $(\Cal H,\delta)$ and $\Psi$ correspond to each other as in
1.3 and 1.4 then $(\pi_0,b)$ and $\psi$ correspond to each other via the GNS construction.

\vskip .1in
\noindent
{\bf 1.6. Examples from free probability}.
From above we have two main examples of closable derivations, those which are inner, and those which 
come from cocycles on groups.  In [V1] and [V2] Voiculescu uses certain derivations in a key role for
his non-microstates approach to free entropy and mutual free information.  We will recall these derivations
which will give us more examples of closable derivations under certain circumstances.

\vskip .1in
\noindent
{\it 1.6.1 The derivation $\partial_X$ from} [V1].
Suppose $B \subset N$ is a $*$-subalgebra with $1 \in B$ and $X = X^* \in N$.  If we denote by $B[X]$ the subalgebra
generated by $B$ and $X$, and if $X$ and $B$ are algebraically free (i.e. they do not satisfy any nontrivial algebraic relations) then there is a well-defined unique
derivation 
$$
\partial_X : B[X] \rightarrow B[X] \otimes B[X] \subset L^2(N, \tau) \otimes L^2(N, \tau)
$$
such that $\partial_X(X) = 1 \otimes 1$ and $\partial_X(b) = 0$ $\forall b \in B$.

We note that if $\partial_X$ is inner then by identifying $L^2(N, \tau) \otimes L^2(N, \tau)$
with the Hilbert-Schmidt operators we would have that there exists a Hilbert-Schmidt operator which
commutes with $B$.  Therefore if $B$ contains a diffuse element 
(i.e. an element which generates a von Neumann algebra without minimal projections) then we must have that
$\partial_X$ is not inner.

Recall from [V1] that the conjugate variable $J(X:B)$ of $X$ w.r.t. $B$ is an element in 
$L^1(W^*(B[X]), \tau)$ such that 
$\tau( J(X:B) m) = \tau \otimes \tau (\partial_X(m))$ $\forall m \in B[X]$,
i.e. $J(X:B) = \partial_X^*(1 \otimes 1)$.

If $J(X:B)$ exists and is in $L^2(N, \tau)$ 
(as in the case when we perturb a set of generators by
free semicircular elements) then by Corollary 4.2 in [V1] we have that  
$\partial_X$ is a closable derivation.

\vskip .1in
\noindent
{\it 1.6.2 The derivation $\delta_{A:B}$ from} [V2].
Suppose $A, B \subset N$ are two $*$-subalgebras with $1 \in A,B$.  
If we denote by $A \vee B$ the subalgebra generated by 
$A$ and $B$, and if $A$ and $B$ are algebraically free then we may define a unique derivation
$$
\delta_{A:B} : A \vee B \rightarrow (A \vee B) \otimes (A \vee B) \subset L^2(N, \tau) \otimes L^2(N, \tau)
$$
by $\delta_{A:B}(a) = a \otimes 1 - 1 \otimes a$ $\forall a \in A$, and
$\delta_{A:B}(b) = 0$ $\forall b \in B$.

We note that for the same reason as above if $B$ contains a diffuse element and $A \not= \Bbb C$ then we must 
have that the derivation is not inner.

Recall from [V2] that the liberation gradient $j(A:B)$ of $(A,B)$ is an element in
$L^1(W^*(A \cup B), \tau)$ such that 
$\tau( j(A:B) m ) = \tau \otimes \tau (\delta_{A:B}(m))$ $\forall m \in A \vee B$,
i.e. $j(A:B) = \delta_{A:B}^*(1 \otimes 1)$.

If $j(A:B)$ exists and is in $L^2(N, \tau)$ then by Corollary 6.3 in [V2]
$\delta_{A:B}$ is a closable derivation.

\vskip .1in
\noindent
{\bf 1.7. Generators of completely positive semigroups}.
Let $N$ be a finite von Neumann algebra with normal faithful trace $\tau$.
A weak*-continuous semigroup $\{\phi_t\}_{t \geq 0}$ on $N$ is said to be
symmetric if $\tau(x\phi_t(y)) = \tau(\phi_t(x)y)$, $\forall x,y \in N$, and
completely Markovian if each $\phi_t$ is a unital c.p. map on $N$.
We denote by $\Delta$ the generator of a symmetric completely Markovian semigroup 
$\{\phi_t\}_{t \geq 0}$ on $N$, i.e. $\Delta$ is the 
densely defined operator on $N$ described by
$D(\Delta) = \{x \in N : {x - \phi_t(x) \over t} \ {\text { has a weak limit as }} \ t \rightarrow 0 \}$,
and $\Delta(x) = \lim_{t \rightarrow 0} {x - \phi_t(x) \over t}$, we also denote by $\Delta$ the generator
of the corresponding semigroup on $L^2(N, \tau)$.
Then $\Delta$ describes a completely Dirichlet form [DL] on $L^2(N, \tau)$ by
$D(\Cal E) = D(\Delta^{1/2})$, $\Cal E(x) = \| \Delta^{1/2}(x) \|_2^2$.

From [DL] we have that $D(\Cal E) \cap N$ is a weakly dense $*$-subalgebra and hence it follows from 
[S1] that there exists a Hilbert $N$-$N$ bimodule $\Cal H$ and a closeable derivation
$\delta: D(\Cal E) \cap N \rightarrow \Cal H$ such that 
$\Cal E(x) = \| \delta(x) \|^2$, $\forall x \in D(\Cal E) \cap N$.
Conversely it follows from [S2] that if $D(\delta)$ is a weakly dense $*$-subalgebra with $1 \in D(\delta)$
and $\delta: D(\delta) \rightarrow \Cal H$ is a closable derivation, then 
the closure of the quadratic form given by $\| \delta(x) \|^2$ is completely Dirichlet on $L^2(N, \tau)$
and hence generates a symmetric completely Markovian semigroup as above (see also [CiSa]).

From sections 1.3 and 1.4, and from the remarks above we obtain the following.

\proclaim{1.7.1. Theorem} Let $N_0 \subset N$ be a weakly dense $*$-subalgebra with $1 \in N_0$ and suppose
$\Psi: N_0 \rightarrow L^1(N, \tau)$ is a closable, conservative, symmetric c.c.n. map, such that
$\Psi^{-1}(L^2(N, \tau))$ is weakly dense in $N$.  
Then $\Delta = \Psi|_{\Psi^{-1}(L^2(N, \tau))}$ is closable as a densely defined operator on $L^2(N, \tau)$
and $\overline \Delta$ is the generator of a symmetric completely Markovian
semigroup on $N$.
Conversely if $\Delta$ is the generator of a symmetric completely Markovian
semigroup on $N$ then $\Delta$ extends to a conservative, symmetric c.c.n. map
$\Psi: N_0 \rightarrow L^1(N, \tau)$ where $N_0$ is the $*$-subalgebra generated
by $D(\Delta)$.
\endproclaim

\vskip .1in
\head 2. A characterization of inner derivations. \endhead

\vskip .1in
Let $N$ be a finite von Neumann algebra with normal faithful trace $\tau$.
Given a symmetric c.c.n. map $\Psi$ on $N$ we will now give a characterization of when $\Psi$ is norm bounded.

\proclaim{2.1. Theorem} Let $\Psi:D_\Psi \rightarrow L^1(N, \tau)$ 
be a closable, conservative, symmetric c.c.n. map with weakly dense
domain $D_\Psi$.  Let $\delta:D_\Psi \rightarrow \Cal H$ be the closable derivation
associated with $\Psi$.  
Then the following conditions are equivalent:
\vskip .05in
\noindent
(a) $\delta$ extends to an everywhere defined derivation $\delta'$ which is inner and such that
given any $x \in N$ there exists a constant $C_x > 0$ such that 
$|\langle \delta'(x),\delta'(y) \rangle| \leq C_x \|y\|_1$, $\forall y \in N$.
\vskip .05in
\noindent
(b) There exists a constant $C > 0$ such that 
$|\langle \delta(x),\delta(y) \rangle| \leq C\|x\| \|y\|_1$, $\forall x,y \in D_\Psi$.
\vskip .05in
\noindent
(c) $\Psi$ is norm bounded on $(D_\Psi)_1$.
\vskip .05in
\noindent
(d) The image of $\Psi$ is contained in $N \subset L^1(N, \tau)$ and 
$-\Psi$ extends to a mapping which generates a norm continuous semigroup of normal c.p. maps.
\vskip .05in
\noindent
(e) There exists $k \in N$ and a normal c.p. map $\phi:N\rightarrow N$ such that $\Psi(x) = k^*x + xk - \phi(x)$,
$\forall x \in D_\Psi$.
\endproclaim
\vskip .05in
\noindent
{\it Proof}.  
(a) $\Rightarrow$ (c):  
Let $\delta'$ be the everywhere defined extension of $\delta$, and
let $\Psi'$ be the c.c.n. map associated with $\delta'$.  Since given any $x \in N$ 
there exists a constant $C_x > 0$ such that 
$|\langle \delta'(x),\delta'(y) \rangle| \leq C_x \|y\|_1$, $\forall y \in N$ we have that the 
image of $\Psi'$ is contained in $N$.
Also since $\Psi'(1) = 0$ we have that for all $x \in D_{\Psi'}$, 
$\Psi'(x^*x) - x^*\Psi'(x) - \Psi'(x^*)x \leq 0$ and so $-\Psi'$ is a dissipation ([L],[Ki]).
As $-\Psi'$ is also everywhere defined, it is bounded by Theorem 1 in [Ki].
\vskip .05in
\noindent
(b) $\Leftrightarrow$ (c):
Suppose (b) holds, then for all $x,y \in D_\Psi$,
$$
|\tau(\Psi(x)y^*)| = |\langle \delta(x),\delta(y) \rangle| \leq C\|x\| \|y\|_1.
$$
So by taking the supremum over all $y \in D_\Psi$ such that $\|y\|_1 \leq 1$ we have that
$\|\Psi(x)\| \leq C\|x\|$, $\forall x \in D_\Psi$.

Suppose now that $\Psi$ is bounded by $C > 0$.  Then for all $x,y \in D_\Psi$,
$$
|\langle \delta(x),\delta(y) \rangle| = |\tau(\Psi(x)y^*)| \leq \| \Psi(x)\| \|y\|_1 \leq C \|x\| \|y\|_1.
$$
\vskip .05in
\noindent
(c) $\Rightarrow$ (d):  This follows from [E] Proposition 2.10.
\vskip .05in
\noindent
(d) $\Rightarrow$ (e):  This is Theorem 3.1 in [ChE].
\vskip .05in
\noindent
(e) $\Rightarrow$ (a):  Suppose that for $k \in N$ and $\phi$ c.p. we have 
$\Psi(x) = k^*x + xk - \phi(x)$, $ \forall x \in N$.
Let $\phi' = \tau(\phi(1))^{-1}\phi$ and let $(\Cal H, \xi)$ be the pointed Hilbert $N$-bimodule
associated with $\phi'$.  Hence if we set $\delta'(x) = (\tau(\phi(1))/2)^{1/2}[x,\xi]$ then
we have $\delta' \cong \delta$.
By replacing $k$ with ${1 \over 2}(k + k^*)$ and $\phi$ with ${1 \over 2}(\phi + \phi^*)$ we may
assume that $\phi$ is symmetric, it is then easy to verify that there exists a constant $C > 0$
such that for all $x,y \in N$, $|\langle \delta'(x),\delta'(y)\rangle| \leq C\|x\| \|y\|_1$.
Hence $\delta'$ gives an everywhere defined extension of $\delta$ which satisfies the required properties.
\hfill $\square$

\vskip .1in
Our next result is in the same spirit as Theorem 2.1 and provides several equivalent conditions 
for when a closable derivation is inner.

\proclaim{2.2. Theorem} Let $\Psi:D_\Psi \rightarrow L^1(N, \tau)$ 
be a closable, conservative, symmetric c.c.n. map with weakly dense
domain $D_\Psi$.  Let $\delta:D_\Psi \rightarrow \Cal H$ be the closable derivation
associated with $\Psi$. Then the following
conditions are equivalent:
\vskip .05in
\noindent
$(\alpha)$ $\delta$ is inner.
\vskip .05in
\noindent
$(\beta)$ $\delta$ is bounded on $(D_\Psi)_1$.
\vskip .05in
\noindent
$(\gamma)$ $\Psi$ is $\|\cdot\|_1$-bounded on $(D_\Psi)_1$.
\vskip .05in
\noindent
$(\delta)$ $\Psi$ can be approximated uniformly by c.p. maps in the following sense:  
for all $\varepsilon > 0$, there exists $k \in N$, and $\phi$ a normal c.p. map such that 
$\|\Psi(x) - k^*x - xk + \phi(x)\|_1 \leq \varepsilon \|x\|$, $ \forall x \in D_\Psi$.
\endproclaim
\vskip .05in
\noindent
{\it Proof}. 
$(\alpha)$ $\Rightarrow$ $(\delta)$:
Suppose $\xi \in \Cal H$ such that $\delta(x) = x\xi - \xi x$, $ \forall x \in D_\Psi$.
Let $\varepsilon > 0$.
Since the subspace of ``left and right bounded'' vectors is dense in $\Cal H$,
let $\xi_0 \in \Cal H$ such that there exists a constant $C >0$ such that 
$\|x\xi_0\| \leq C\|x\|_2$, $\forall x \in N$, $\| \xi_0 \| \leq \| \xi \|$, and also $\|\xi - \xi_0\| < \varepsilon/8\|\xi\|$.
As $\xi_0$ is ``bounded'' we may let $\phi_{\xi_0}$ be the normal c.p. map associated with
$\xi_0/\|\xi_0\|$.  Let $\phi = 2 \|\xi_0\|^2 \phi_{\xi_0}$, and let $k = \phi(1)/2$.

Note that since $\delta$ is real we have that $\xi_0$ is also real, i.e. 
$\langle x\xi_0,\xi_0 y \rangle = \langle y^*\xi_0,\xi_0 x^* \rangle$,
$\forall x,y \in N$.

Then if $x, y \in D_\Psi$ we have:
$$
\tau( (\Psi(x) - k^*x - xk + \phi(x))y^*)
$$
$$
= \tau(\Psi(x)y^*) - {1 \over 2}\tau(\phi(1)xy^*) - {1 \over 2}\tau(x\phi(1)y^*) + \tau(\phi(x)y^*)
$$
$$
= \langle \delta(x), \delta(y) \rangle - \langle xy^*\xi_0,\xi_0 \rangle -
  \langle y^*x\xi_0,\xi_0 \rangle + 2\langle x\xi_0 y^*,\xi_0 \rangle
$$
$$
= \langle x\xi - \xi x, y\xi - \xi y \rangle - \langle x\xi_0 - \xi_0 x, y\xi_0 - \xi_0 y \rangle.
$$
Hence:
$$
|\langle \Psi(x) - k^*x - xk + \phi(x),y \rangle|
$$
$$
\leq \|x\xi - \xi x\| \|y\xi - \xi y - y\xi_0 + \xi_0 y\| + \|y\xi_0 - \xi_0 y\| \|x\xi - \xi x - x\xi_0 + \xi_0 x\|
$$
$$
\leq 4 \| x \| \| \xi \| \| y \| \| \xi - \xi_0 \| + 4 \| y \| \| \xi_0 \| \| x \| \| \xi - \xi_0 \|
$$
$$
\leq \varepsilon \|x\| \|y\|.
$$
Thus by taking the supremum over all $y \in (D_\Psi)_1$ we have the desired result.
\vskip .05in
\noindent
$(\delta)$ $\Rightarrow$ $(\gamma)$:
Let $k \in N$ and $\phi$ c.p. such that $\|\Psi(x) - k^*x - xk + \phi(x)\|_1 \leq \|x\|$.
By ([P2] 1.1.2) $\|\phi(x)\|_2 \leq \|\phi(1)\|_2 \|x\|$, $\forall x \in N$. 
Hence for all $x \in D_\Psi$:
$$
\|\Psi(x)\|_1 \leq \|\Psi(x) - k^*x - xk + \phi(x)\|_1 + \|k^*x - xk + \phi(x)\|_2
$$
$$
\leq (1 + 2\|k\|_2 + \|\phi(1)\|_2)\|x\|.
$$
Thus $\Psi$ is bounded in $\|\cdot\|_1$ on $(D_\Psi)_1$.
\vskip .05in
\noindent
$(\gamma)$ $\Rightarrow$ $(\beta)$:
Suppose $\|\Psi(x)\|_1 \leq C\|x\|$, $\forall x \in D_\Psi$, then:
$$
\|\delta(x)\|^2 = \tau(\Psi(x)x^*) \leq \|\Psi(x)\|_1 \|x\| \leq C\|x\|^2,
$$
for all $x \in D_\Psi$.
\vskip .05in
\noindent
$(\beta)$ $\Rightarrow$ $(\alpha)$:
As $\delta$ is bounded on $(D_\Psi)_1$ we may extend $\delta$ to a derivation on the 
$C^*$-algebra $A$ which is generated by $D_\Psi$.  
Let $X = \{\delta(u)u^* | u \in \Cal U(A)\}$, for each 
$v \in \Cal U(A)$ we let $v$ act on $\Cal H$ by $v\cdot \xi = v\xi + \delta(v)$. 
Let $\xi_0$ be the center of the set $X$.  Then since 
$\forall \xi,\eta \in \Cal H$, 
$\|v\cdot \xi - v\cdot \eta \| = \|\xi - \eta \| = \|\xi v - \eta v\|$, 
the center of the set $v\cdot X$ is $v\cdot \xi_0$, and
the center of the set $Xv$ is $\xi_0 v$.  Further we have that $v\cdot X = Xv$ and thus
$v\cdot \xi_0 = \xi_0 v$.  Since $v$ was arbitrary and every $x \in A$ is a linear combination
of unitaries we have that $\delta(x) = \xi_0 x - x\xi_0$, $ \forall x \in A$.
\hfill $\square$

\vskip .1in
In general we may have that $\Psi$ is unbounded in $\|\cdot\|_1$ even if $\|\phi_t(x) - x\|_2$ converges to
0 uniformly on $N_1$.  However we will show in the next section that if $N$ has property (T) and
the domain of $\Psi$ contains a ``critical set'' as in 
Proposition 1 of [CJ] then this cannot happen.  

\vskip .1in
\head 3. Property (T) in terms of closable derivations. \endhead

\vskip .1in
Given a finite von Neumann algebra $M$ with countable decomposable center.  We will say that
$M$ has property (T) if the inclusion $(M \subset M)$ is rigid in the sense of [P2], i.e.
$M$ has property (T) if and only if there exists a normal faithful tracial state $\tau'$ on $M$
such that one of the following equivalent conditions hold:

\vskip .05in
\noindent
1.  $\forall \varepsilon > 0$, $\exists F'  \subset M$ finite and $\delta' > 0$ such that if 
$\Cal H$ is a Hilbert $M-M$ bimodule with a vector $\xi \in \Cal H$ satisfying the conditions 
$\| \langle \cdot \xi, \xi \rangle - \tau' \| \leq \delta'$,
$\| \langle \xi \cdot, \xi \rangle - \tau' \| \leq \delta'$,
and $\| y \xi - \xi y \| \leq \delta'$, $\forall y \in F'$
then $\exists \xi_0 \in \Cal H$ such that $\| \xi_0 - \xi \| \leq \varepsilon$ 
and $x \xi_0 = \xi_0 x$, $\forall x \in M$.

\vskip .05in
\noindent
2.  $\forall \varepsilon > 0$, $\exists F \subset M$ finite and $\delta > 0$ such that if
$\phi: M \rightarrow M$ is a normal, completely positive map with 
$\tau' \circ \phi \leq \tau'$, $\phi(1) \leq 1$ and 
$\| \phi(y) - y \|_2 \leq \delta$, $\forall y \in F$, then
$\| \phi(x) - x \|_2 \leq \varepsilon$, $\forall x \in M$, $\|x\| \leq 1$.

\vskip .05in
Furthermore it was shown in [P2] that the above definition is independent of the trace $\tau'$, and in the 
case when $N$ is a factor this agrees with the original definition in [CJ].

In this section we will obtain a characterization of property (T) in terms of certain boundedness
conditions on c.c.n. maps.  As we are dealing with unbounded maps the domain of a map will 
be of crucial importance.  We will thus want to consider c.c.n. maps whose domain contains 
a ``critical set'', which by Remark 4.1.6 of [P1] motivates the following.

\proclaim{3.1. Definition}  Suppose that $N$ is a {\rm II}$\sb 1$ factor and
let $N_0 \subset N$ be a weakly dense $*$-subalgebra of $N$, $1 \in N_0$, 
then $N_0$ contains a non-$\Gamma$ set if $\exists F \subset N_0$ finite, $K > 0$ such that
$\forall \xi \in L^2(N, \tau)$, if $\langle \xi, 1 \rangle = 0$ then we have 
$\| \xi \|_2^2 \leq K \Sigma_{x \in F} \| x \xi - \xi x \|_2^2$.
\endproclaim

\vskip .05in
Note that by Lemma 2.4 of [C1] one can check that $N_0$ has a non-$\Gamma$ set if and only if 
there exists a finitely generated subgroup $\Cal G \subset {\text { Int }} C^*(N_0)$ such that
there is no non-normal $\Cal G$-invariant state on $N$.  Also it follows from the definition
that $N \subset N$ contains a non-$\Gamma$ set if and only if $N$ does not have property $\Gamma$
of Murray and von Neumann [MvN].
Also note that if $\Lambda$ is a countable ICC group then by [Ef]
$\Lambda$ is not inner amenable if and only if 
$\Bbb C\Lambda$ contains a non-$\Gamma$ set.

\vskip .1in
We now come to the main result which is to give several equivalent characterizations of property (T),
in particular we obtain a 1-cohomology characterization of property (T) which is the analogue of 
the Delorme-Guichardet Theorem from group theory.

\proclaim{3.2. Theorem}  Suppose that $N$ is a separable 
finite factor with normal faithful trace $\tau$.  Let $N_0 \subset N$ be a weakly dense $*$-subalgebra
such that $1 \in N_0$ and $N_0$ is countably generated as a vector space.
Consider the following conditions:

\vskip .05in
\noindent
(a) $N$ has property (T).

\vskip .05in
\noindent
(b) $\exists F \subset N_0$ finite and $K > 0$ such that if $\Cal H$ is a Hilbert 
$N$-$N$ bimodule, $\xi \in \Cal H$,
and if $\delta_\xi = \max_{x \in F} \{ \| x\xi - \xi x \| \}$
then $\exists \xi_0 \in \Cal H$ such that 
$x \xi_0 = \xi_0 x, \forall x \in N$ and
$\| \xi_0 - \xi \| \leq \delta_\xi K$.

\vskip .05in
\noindent
(c) Every densely defined closable derivation on $N_0$ is inner.

\vskip .05in
\noindent
(d) Every closable, conservative, symmetric c.c.n. map on $N_0$ is bounded in $\| \cdot \|_1$ on $(N_0)_1$. 

\vskip .05in
\noindent
(e)  $\exists F' \subset N_0$ finite and $K' >0$ such that if $\phi : N \rightarrow N$
is a c.p. map with $\phi(1) \leq 1$, $\tau \circ \phi \leq \tau$, $\phi = \phi^*$, 
and if $\delta_\phi' = \max_{x \in F'} \{ \| x - \phi(x) \|_2 \}$
then $\tau((y - \phi(y))y^*) \leq K' \delta_\phi'$, $\forall y \in (N)_1$.

Then (b) $\Rightarrow$ (c) $\Rightarrow$ (d) $\Rightarrow$ (e) $\Rightarrow$ (a).
If moreover $N_0$ contains a non-$\Gamma$ set
then we also have (a) $\Rightarrow$ (b).
\endproclaim

\vskip .05in
\noindent
{\it Proof}.
(b) $\Rightarrow$ (c).
Let $\delta: N_0 \rightarrow \Cal H$ be a closable derivation, note that by 1.2 we may assume
that $\delta$ is real.
Let $\phi_t : N \rightarrow N$ be the semigroup of normal symmetric c.p. maps associated with $\delta$.
Then $\forall y \in N_0$ 
$\| \delta(y) \|^2 = \lim_{t \rightarrow \infty} \tau({y - \phi_t(y) \over t}y^*)$.

Let $(\Cal H_t, \xi_t)$ be the pointed correspondence obtained from $\phi_t$, then since $\phi_t$ is 
unital and symmetric $\| y\xi_t - \xi_t y \|_2^2 = 2\tau((y - \phi_t(y))y^*)$ $\forall y \in N$.
Let $F \subset N_0$ and $K > 0$ be as in (b) and let 
$C = \sup_{0 < t \leq 1, x \in F} \tau((x - \phi_t(x))x^*)/t$.
Then:
$$
\| \delta(y) \|^2 = \lim_{t \rightarrow 0} \tau({y - \phi_t(y) \over t}y^*)
$$
$$
= \lim_{t \rightarrow 0} \|y \xi_t - \xi_t y\|_2^2/2t
$$
$$
\leq  2\sup_{0 < t \leq 1, x \in F} \|x \xi_t - \xi_t x\|_2^2 K^2/t 
$$
$$
= 4 \sup_{0 < t \leq 1, x \in F} \tau((x - \phi_t(x))x^*) K^2/t = 4CK^2,
$$
for all $y \in (N_0)_1$.  Thus by Theorem 2.2 
we have that $\delta$ is inner.

\vskip .05in
\noindent
(c) $\Rightarrow$ (d).  
This follows from Theorem 2.2.

\vskip .05in
\noindent
(d) $\Rightarrow$ (e).
Let $\{x_n\}_{n \in \Bbb N}$ be a sequence in $(N_0)_1$ such that
$N_0 = \text{sp}\{x_n\}_{n \in \Bbb N}$.  If (e) does not hold then for each $k \in \Bbb N$ there
exists a c.p. map $\phi_k: N \rightarrow N$ such that $\phi_k(1) \leq 1$, $\tau \circ \phi_k \leq \tau$,
$\phi_k = \phi_k^*$, 
and there exists $y_k \in (N_0)_1$ such that
$\tau((y_k - \phi_k(y_k))y_k^*) > 4^k \delta_k'$ where
$\delta_k' = \max_{j \leq k} \{ \| x_j - \phi_k(x_j) \|_2 \}$.

Let $\Psi_k = (\text{id} - \phi_k)/ \delta_k'$, and let 
$\Psi = \Sigma_{k=1}^\infty 2^{-k} \Psi_k$.
Then since $N_0 = \text{sp}\{x_n\}_{n \in \Bbb N}$, $\Psi: N_0 \rightarrow L^2(N, \tau)$ is a well 
defined symmetric c.c.n. map with $\Psi(1) \geq 0$.
Also since $\phi_k(1) \leq 1$, $\tau \circ \phi_k \leq \tau$, and $\phi_k = \phi_k^*$ 
if we let $(\Cal H_k, \xi_k)$ be the pointed Hilbert $N$-$N$ bimodule
corresponding to $\phi_k$ then we have
$2\tau( (x - \phi_k(x)) x^*) \geq \tau \circ \phi_k (x^*x) + \tau(x^*x \phi_k(1)) - 2 \tau(\phi(x)x^*)
				= \| x \xi_k - \xi_k x \|^2 \geq 0$, 
$\forall x \in N$.
Thus
$$
\| \Psi(y_k) \|_1 \geq \tau(\Psi(y_k)y_k^*) 
$$
$$
\geq 2^{-k} \tau((y_k - \phi_k(y_k)y_k^*)/\delta_k' > 2^{k},
$$
for all $k \in \Bbb N$.
Hence if we let
$\Psi'(x) = \Psi(x) - x\Psi(1)/2 - \Psi(1)x/2$ then
$\Psi'$ is a closable, conservative, symmetric c.c.n. map which is unbounded
in $\|\cdot\|_1$ on $(N_0)_1$.

\vskip .05in
\noindent
(e) $\Rightarrow$ (a).  
Let $F'$ and $K'$ be as in (e) and let $\varepsilon > 0$.
Suppose $\phi: N \rightarrow N$ is a c.p. map such that 
$\phi(1) \leq 1$, $\tau \circ \phi \leq \tau$, $\phi = \phi^*$,
and $\| x - \phi(x) \|_2 < \varepsilon^2/2K'$, $\forall x \in F'$.
Let $(\Cal H_\phi, \xi_\phi)$ be the pointed Hilbert $N$-$N$ bimodule 
associated with $\phi$, then since $\| \phi(1) \|_2 \leq 1$ by Lemma 1.1.3 of [P2]
we have that 
$\| y - \phi(y) \|_2^2 \leq  \| y \xi_\phi - \xi_\phi y \|^2
			= 2  \tau((y - \phi(y))y^*)
			\leq 2K'\delta_\phi' < \varepsilon^2$,
$\forall y \in (N)_1$.
Hence by Lemma 3 of [PeP] $N$ has property (T).

\vskip .05in
\noindent
(a) $\Rightarrow$ (b).
If $N_0$ contains a non-$\Gamma$ set then 
Remark 4.1.6 of [P1] shows that 
one can then apply Proposition 1 of [CJ]
to obtain the desired result.
\hfill $\square$

\vskip .1in

We note that (d) $\Rightarrow$ (e) in Theorem 3.2 can be suitably adapted to the case of inclusions of 
$\sigma$-compact and locally compact groups thus showing that an inclusion of groups has 
relative property (T) if and only if ``$\delta$ depends linearly on $\varepsilon$'', answering a
question of Jolissaint (see Theorem 1.2 in [J]).

\vskip .1in
Let $B \subset N$, $1 \in B$ be a $*$-subalgebra and $X = X \in N$. 
Recall from [V1] that a dual operator to $(X;B)$ in $L^2(N, \tau)$ is an 
operator $Y \in \Cal B(L^2(N, \tau))$ such that
$$
[B, Y] = 0 \ {\text{  and  }} \ [X, Y] = P_1 
$$ 
where $P_1$ is the orthogonal projection onto $\Bbb C1$.

\proclaim{3.3. Corollary}  Suppose that $N$ is a separable finite factor with property (T), let 
$B \subset M$, $1 \in B$ be a $*$-subalgebra and $X = X^* \in N$ such that 
$B[X]$ generates $N$ as a von Neumann algebra.  Suppose that $B$ contains a diffuse element 
and $B[X]$ contains a non-$\Gamma$ set, then 
the conjugate variable $J(X:B)$ does not exist
in $L^2(N, \tau)$, i.e. $\Phi^*(X:B) = \infty$.
Also $(X;B)$ does not have a dual operator in $L^2(N, \tau)$.
\endproclaim

\vskip .05in
\noindent
{\it Proof.} If the conjugate variables $J(X:B)$
did exist in $L^2(N, \tau)$ then as in 1.6.1 we would have a closeable derivation on $B[X]$
which is not inner.  Therefore by Theorem 3.2 this cannot happen.

The fact that $(X;B)$ does not have a dual system in $L^2(N, \tau)$ then follows directly from [V1].
\hfill $\square$

\vskip .1in
\head 4. Property (T) and amalgamated free products. \endhead

\vskip .1in
We include here an application of the above ideas, showing that a large class of 
amalgamated free products do not have property (T).  
We first prove that if $N$ has property (T) then even though a c.c.n. maps may be unbounded on some domains it must still satisfy a certain 
condition on its rate of growth.

\proclaim{4.1. Theorem} Suppose $N$ is a finite von Neumann algebra with normal faithful tracial state $\tau$,
if $N$ has property (T) and 
$\Psi: D_\Psi \rightarrow L^2(N,\tau) \subset L^1(N, \tau)$ 
is a conservative, symmetric c.c.n. map, then
given any sequence $\{ x_n \}_n$ in $(D_\Psi)_1$ such that 
$\| \Psi(x_n) \|_2 \rightarrow \infty$, we have that $\| \Psi(x_n) \|_2 / \| \Psi(x_n) \| \rightarrow 0$.
\endproclaim 

\vskip .05in
\noindent
{\it Proof.}
Let $\{ \Phi_t \}_t$ be the semigroup of unital normal symmetric c.p. maps associated with $\Psi$ as in 1.7, and for each $\beta > 0$ let $\varepsilon_\beta = \sup_{t \leq \beta, x \in N_1} \| \Phi_t(x) - x \|_2$.  Since $N$ has property (T) we have that $\varepsilon_\beta \rightarrow 0$ as $\beta \rightarrow 0$.

For all $\beta > 0$, and $x \in (D_\Psi)_1$ we have:

$$
\int\limits_0^{\beta} \Phi_t \circ \Psi (x) dt 
= \lim_{s \rightarrow 0} \int\limits_0^{\beta} \Phi_t ((\Phi_s(x) - x)/ s) dt
$$ 
$$
= \lim_{s \rightarrow 0} {1 \over s} ( \int\limits_0^{\beta} \Phi_{t+s}(x) dt - 
\int\limits_0^\beta  \Phi_t(x) dt)
$$
$$
= \lim_{s \rightarrow 0} {1 \over s} (\int\limits_s^{\beta + s} \Phi_t (x) dt - 
\int\limits_0^\beta  \Phi_t(x) dt)
$$
$$
= \lim_{s \rightarrow 0} {1 \over s} (\int\limits_\beta^{\beta + s} \Phi_t (x) dt - 
\int\limits_0^s  \Phi_t(x) dt)
$$
$$
= \Phi_\beta (x) - x.
$$

Hence for all $x \in (D_\Psi)_1$:

$$
\| \Psi(x) \|_2 \leq \| {1 \over \beta} \int\limits_0^{\beta} \Phi_t \circ \Psi (x) dt \|_2 +
     \| {1 \over \beta} \int\limits_0^\beta (\Phi_t \circ \Psi(x) - \Psi(x)) dt \|_2
$$
$$
\leq {\varepsilon_\beta \over \beta} 
+  {1 \over \beta} \int\limits_0^\beta \| (\Phi_t \circ \Psi(x) - \Psi(x)) \|_2 dt
$$
$$
\leq {\varepsilon_\beta \over \beta} 
+ \| \Psi(x) \| {1 \over \beta} \int\limits_0^\beta \varepsilon_t dt
$$
$$
\leq {\varepsilon_\beta \over \beta} + \| \Psi(x) \| \varepsilon_\beta.
$$

Thus $\varepsilon_\beta \geq \| \Psi(x) \|_2 \beta / (1 + \| \Psi(x) \| \beta)$ and since 
$\varepsilon_\beta \rightarrow 0$ the result follows.
\hfill $\square$
\vskip .1in

\proclaim{4.2. Corollary} Let $N_1$ and $N_2$ be finite von Neumann algebras with normal faithful 
tracial states $\tau_1$ and $\tau_2$ respectively, suppose that $B$ is a common
von Neumann subalgebra such that $\tau_1|_B = \tau_2|_B$.  Suppose also that there are unitaries
$u_i \in \Cal U (N_i)$ such that $E_B(u_i) = 0$, $i = 1,2$.
Then $M = N_1 *_B N_2$ does not have property (T).
\endproclaim

\vskip .05in
\noindent
{\it Proof.}
Let $\tau = \tau_1 *_B \tau_2$ be the trace for $M$ and let $\Cal H = L^2(M, \tau) \otimes_B  L^2(M, \tau)$.  Define 
$\delta$ to be the unique derivation from the algebraic amalgamated free product to $\Cal H$ 
which satisfies $\delta(a) = a \otimes_B 1 - 1  \otimes_B a$, $\forall a \in N_1$, and
$\delta(b) = 0$, $\forall b \in N_2$.
By Corollary 5.4 in [NShSp] $\delta^*(1 \otimes_B 1) = 0$ and so in particular just as in the non-amalgamated case
we have that $\delta$ is a closable derivation and furthermore if $u_1, u_2$ are the unitaries as above,
and $z \in N_0$ then:

$$
\langle  \delta((u_1u_2)^n ), \delta(z) \rangle 
= \Sigma_{j = 0}^{n-1} \langle (u_1u_2)^j u_1 \otimes_B u_2(u_1u_2)^{n-j-1} 
                                                   - (u_1u_2)^j \otimes_B (u_1u_2)^{n-j},
                                                        \delta(z) \rangle
$$
$$
= \Sigma_{j=0}^{n-1} \langle 1 \otimes_B 1, u_1^* (u_2^*u_1^*)^j \delta(z) (u_2^*u_1^*)^{n-j-1} u_2^*
                                                                           - (u_2^*u_1^*)^j \delta(z) (u_2^*u_1^*)^{n-j} \rangle,
$$
also for each $0 \leq j < n$, by using the Leibnitz rule for the derivation we may rewrite 
$(u_2^*u_1^*)^j \delta(z) (u_2^*u_1^*)^{n-j}$ as a sum of three terms: 

$$
 \delta ((u_2^*u_1^*)^j z (u_2^*u_1^*)^{n-j}) \tag 1
$$
$$
- \Sigma_{k = 0}^{j-1} (u_2^*u_1^*)^k u_2^* \delta(u_1^*) (u_2^*u_1^*)^{j - k - 1} z (u_2^*u_1^*)^{n-j} \tag 2
$$
$$
 - \Sigma_{i = 0}^{n - j-1} (u_2^*u_1^*)^j z (u_2^*u_1^*)^i u_2^* \delta(u_1^*) (u_2^*u_1^*)^{n-j-i-1}, \tag 3
$$
however when we take the inner product with $1 \otimes_B 1$ the first term will be 0 as mentioned above,
and by freeness (since $u_1$, and $u_2$ have expectation 0) the other terms will be 0 except when $i = n - j - 1$ in the third term where we have

$$
- \langle 1 \otimes_B 1,  (u_2^*u_1^*)^{j} z (u_2^*u_1^*)^{n - j - 1} u_2^* \delta(u_1^*) \rangle
= - \langle 1 \otimes_B 1,  (u_2^*u_1^*)^{j} z (u_2^*u_1^*)^{n - j} \otimes_B 1 \rangle
$$
$$
= -\tau(E_B((u_1u_2)^{n-j} z^* (u_1u_2)^j)).
$$
Similarly 
$$
\langle 1 \otimes_B 1, u_1^* (u_2^*u_1^*)^j \delta(z) (u_2^*u_1^*)^{n-j-1} u_2^* \rangle
 = \tau(E_B(u_2(u_1u_2)^{n-j-1}  z^* (u_1u_2)^j u_1)).
$$

Hence from the above equalities we have:
$$
\langle  \delta((u_1u_2)^n ), \delta(z) \rangle 
= \Sigma_{j=0}^{n-1} \tau( u_2(u_1u_2)^{n-j-1}  z^* (u_1u_2)^j u_1 + (u_1u_2)^{n-j} z^* (u_1u_2)^j )
$$
$$
= 2n\tau((u_1u_2)^n z^*).
$$

In particular if $\Psi$ is the c.c.n. map associated with $\delta$ then we have that 
$\Psi((u_1u_2)^n) = 2n(u_1u_2)^n$ and so $\| \Psi ((u_1u_2)^n) \|_2 \rightarrow \infty$ but
$\| \Psi ((u_1u_2)^n) \|_2 / \| \Psi ((u_1u_2)^n) \| \not\rightarrow 0$,
hence by Theorem 4.1 $M$ does not have property (T).
\hfill $\square$

\vskip .1in
Note that the only place where we used the fact that $u_1$ and $u_2$ were unitairies was to 
insure that $2n \| (u_1u_2)^n \|_2 \rightarrow \infty$.  Also note that the conditions of the above Corollary
are satisfied when $M$ is a free product (with amalgamation over $\Bbb C$) as well as when $M$ is a group 
von Neumann algebra coming from an amalgamated free products of groups.
We also mention that from the calculation above we are able to compute explicitly the semigroup of c.p. maps 
that $\delta$ generates, it is the semigroup given by 
$\phi_t = (e^{-2t}{\text {id}} + (1 - e^{-2t})E_B) *_B {\text {id}}$. 

\vskip .1in
\noindent
{\bf Acknowledgments}.  The author would like to thank Professor Sorin Popa for his encouragement and many stimulating discussions.  The author would also like to thank Professor Dima Shlyakhtenko for showing him the references [S1] and [S2], and Professor Jean-Luc Sauvageot for his useful comments on an earlier version of this paper.

\vskip .1in
\head  References\endhead

\item{[AW]} C.A. Akemann, M.E. Walter: {\it Unbounded negative definite functions},
Can, J. Math. {\bf 33} (1981), 862-871.

\item{[BdHV]} B. Bekka, P. de la Harpe, A. Valette: {\it Kazhdan's property
(T)}, Book Preprint, (July 2002).

\item{[ChE]} E. Christensen, D.E. Evans: {\it Cohomology of operator algebras and quantum
dynamical semigroups}, J. London Math. Soc. {\bf 20} (1979), 358-368.

\item{[CiSa]} F. Cipriani, J.-L. Sauvageot: {\it Derivations as square roots of Dirichlet forms},
J. Funct. Anal. {\bf 201} (2003), no. 1, 78-120.

\item{[C1]} A. Connes: {\it Classification of injective factors}, Ann. of Math. {\bf 104}
(1976), 73-115.

\item{[C2]} A. Connes: {\it Classification des facteurs}, Proc. Symp. Pure Math.
{\bf 38} (Amer. Math. Soc., 1982), 43-109.

\item{[C3]} A. Connes: {\it Correspondences}, hand-written notes, 1980.

\item{[CJ]} A. Connes, V.F.R. Jones: {\it Property (T) for von Neumann algebras},
Bull. London Math. Soc. {\bf 17} (1985), 57-62.

\item{[DL]} E. B. Davies, J. M. Lindsay: {\it Non-commutative symmetric Markov semigroups},
Math. Z. {\bf 210} (1992), 379-411.

\item{[De]} P. Delorme: 1-{\it cohomologie des repr\'esentations unitaires de groupes de Lie 
semisimples et r\'esolubles.  Produits tensoriels continus et repr\'esentations},
Bull. Soc. Math. France, {\bf 105} (1977), 281-336.

\item{[Ef]} E. G. Effros: {\it Property $\Gamma$ and inner amenability},
Proc. Amer. Math. Soc. {\bf 47} (1975), 483-486.

\item{[E]} D.E. Evans: {\it Conditionally
completely positive maps on operator algebras}, Quart. J. Math.
Oxford {\bf 4} (1977), 271-284.

\item{[G]} A. Guichardet: {\it Etude de la }1-{\it cohomologie et de la 
topologie du dual pour les groupes de Lie \`a radical ab\'elien},
Math. Ann., {\bf 228} (1977), 215-232.

\item{[J]}  P. Jolissaint: {\it On Property (T) for Pairs of Topological Groups},
preprint 2004.

\item{[Ka]} D.A. Kazhdan: {\it Connection of the dual space of a group with the 
structure of its closed subgroups}, Funct. Anal. Appl. {\bf 1} (1967), 63-65.

\item{[Ki]} A. Kishimoto: {\it Dissipations and derivations},
Commun. math. Phys. {\bf 47} (1976) 25-32.

\item{[L]} G. Lindblad: {\it On the generators of quantum dynamical semigroups},
Commun. math. Phys. {\bf 48} (1976), 119-130.

\item{[MvN]} F. J. Murray, J. von Neumann: {\it On rings of operators {\rm IV}},
Ann. of Math. (2) {\bf 44} (1943), 716-808.

\item{[NShSp]} A. Nica, D. Shlyakhtenko, R. Speicher:  {\it Operator valued distributions. I. 
Characterizations of freeness},  Int. Math. Res. Not. (2002), 1509-1538.

\item{[PeP]} J. Peterson, S. Popa: {\it On the notion of relative property (T) for inclusions of von Neumann algebras},
J. Funct. Anal., {\bf 219} (2005), no. 2, 469-483.

\item{[P1]} S. Popa: {\it Correspondences}, INCREST preprint 1986, unpublished.

\item{[P2]} S. Popa: {\it On a class of type {\rm II}$\sb 1$ factors with Betti numbers invariants},
preprint math.OA/0209130, to appear in Ann. of Math.

\item{[S1]} J.-L. Sauvageot: {\it Tangent bimodules and locality for dissipative operators
on C*-algebras}, Quantum probability and appl., {\rm IV}, Lecture notes in Math. {\bf 1396}
(1989), 322-338. 

\item{[S2]} J.-L. Sauvageot: {\it Quantum Dirichlet forms, differential calculus and semigroups},
Quantum probability and appl., {\rm V}, Lecture notes in Math. {\bf 1442} (1990), 334-346.

\item{[V1]} D. Voiculescu: {\it The analogues of entropy and of Fisher's information measure in free 
probability theory {\rm V}.  Noncommutative Hilbert transform}, Invent. math. {\bf 132} (1998) 189-227.

\item{[V2]} D. Voiculescu: {\it The analogues of entropy and of Fisher's information measure in free 
probability theory {\rm VI}.  Liberation and mutual free information}, Adv. Math. {\bf 146} (1999) 101-166.
\enddocument